\documentclass[12pt,draft]{article}

\usepackage[cp1251]{inputenc}
\usepackage[T2A]{fontenc}
\usepackage[english,russian]{babel}
\usepackage{amsmath,amsfonts,amssymb}
\usepackage{geometry}
\begin{document}

\begin{flushleft}



\textbf{Anna V. Anop and Аleksandr А. Murach}\\\small(Institute of Mathematics of NAS of Ukraine,
Kiev)

\medskip

\large\textbf{REGULAR ELLIPTIC BOUNDARY-VALUE PROBLEMS\\IN THE EXTENDED SOBOLEV SCALE}

\bigskip

\normalsize\textbf{Анна В. Аноп и Александр А. Мурач }\\\small(Институт математики НАН Украины,
Киев)

\medskip

\large\textbf{РЕГУЛЯРНЫЕ ЭЛЛИПТИЧЕСКИЕ КРАЕВЫЕ ЗАДАЧИ\\В РАСШИРЕННОЙ СОБОЛЕВСКОЙ
ШКАЛЕ}\footnote{Исследование поддержано грантом №~01-01-12/2 НАН Украины (в рамках
совместного украинско--российского проекта НАН Украины и Российского фонда
фундаментальных исследований).}

\end{flushleft}

\normalsize

\medskip

We investigate an arbitrary regular elliptic boundary--value problem given in a
bounded Euclidean $C^{\infty}$-domain. We prove that the operator of the problem is
bounded and Fredholm in appropriate pairs of H\"ormander inner product spaces. They
are parametrized with the help of an arbitrary radial function RO-varying at
$\infty$ and form the extended Sobolev scale. We establish a priori estimates for
solutions to the problem and investigate their local regularity on this scale. We
find new sufficient conditions for generalized partial derivatives of the solutions
to be continuous.

\medskip

Досліджено довільну регулярну еліптичну крайову задачу, задану в обмеженій
евклідовій області класу $C^{\infty}$. Доведено, що оператор цієї задачі, є
обмеженим і фредгольмовим у відповідних парах гільбертових просторів Хермандера.
Вони параметризовані за допомогою довільної радіальної функції, RO-змінної на
$+\infty$, та утворюють розширену соболєвську шкалу. Встановлено апріорні оцінки
розв'язків задачі та досліджено їх локальну регулярність в цій шкалі. Знайдені нові
достатні умови неперервності узагальнених частинних похідних розв'язку.

\medskip

\textbf{1. Введение.} В теории эллиптических дифференциальных уравнений важная роль
принадлежит пространствам Соболева. Эллиптические краевые задачи (ЭКЗ) имеют
фундаментальные свойства в соболевских шкалах: фредгольмовость (т.~е. конечность
индекса задачи), априорные оценки решений, локальное повышение регулярности решений
и другие (см., например, монографии [1~-- 7] и обзор \cite{Agranovich97}). С точки
зрения применений этих свойств, в частности, в спектральной теории дифференциальных
операторов, наиболее полезна гильбертова шкала пространств Соболева. Однако, она не
является достаточно тонко градуированной для ряда задач [2, 9~-- 14].

Недавно В.~А.~Михайлецом и вторым автором этой статьи была построена теория
эллиптических операторов  и эллиптических краевых задач для уточненной соболевской
шкалы \cite{MikhailetsMurach10}. Она образована гильбертовыми пространствами
Хермандера $H^{s,\varphi}:=\mathcal{B}_{2,\mu}$ \cite{Hermander63, Hermander83}, для
которых показателем гладкости служит функция
$\mu(\xi):=\langle\xi\rangle^{s}\varphi(\langle\xi\rangle)$ аргумента
$\xi\in\mathbb{R}^{n}$. Здесь $\langle\xi\rangle:=(1+|\xi|^{2})^{1/2}$,
$s\in\mathbb{R}$, а $\varphi$ --- произвольная положительная функция, медленно
меняющаяся на $\infty$ по  Й.~Карамата \cite{Seneta76, BinghamGoldieTeugels89}.
Пространство $H^{s,\varphi}$ на $\mathbb{R}^{n}$ состоит со всех медленно растущих
распределений $w\in\mathcal{S}'(\mathbb{R}^{n})$ таких, что $\mu\widehat{w}\in
L_{2}(\mathbb{R}^{n})$, и определяется на евклидовых  областях и гладких компактных
многообразиях стандартным образом ($\widehat{w}$~--- преобразование Фурье
распределения $w$).

Эта шкала  $\{H^{s,\varphi}\}$ содержит пространства Соболева: $H^{s,1}= H^{s}$,
привязана к соболевской шкале с помощью числового параметра $s$ и тоньше
градуирована, чем последняя. Для применений важно, что каждое пространство
$H^{s,\varphi}$ получается в результате интерполяции с некоторым функциональным
параметром пары соболевских пространств $H^{s_{0}}$ и $H^{s_{1}}$, где
$s_{0}<s<s_{1}$. Поскольку при интерполяции пространств наследуется ограниченность
линейных операторов и их фредгольмовость, то уточненная соболевская шкала оказалась
удобным инструментом в исследовании свойств операторов, порожденных ЭКЗ (см. также
статьи [17~-- 22] и обзор \cite{12BJMA2}).

В этой связи естественно возникает вопрос об описании и возможных применениях класса
всех гильбертовых пространств, интерполяционных относительно пар гильбертовых
пространств Соболева. Из теоремы В.~И.~Овчинникова \cite{Ovchinnikov84} (п.~11.4)
вытекает, что это класс состоит (с точностью до эквивалентности норм) из пространств
Хермандера $H^{\varphi}:=\mathcal{B}_{2,\mu}$, где
$\mu(\xi):=\varphi(\langle\xi\rangle)$, а $\varphi$ --- положительная функция,
RO-меняющаяся на $+\infty$ по В.~Г.~Авакумовичу \cite{Seneta76,
BinghamGoldieTeugels89}. Указанный  класс пространств является расширенной
соболевской шкалой (с помощью интерполяции) \cite{13UMJ3}. Он содержит уточненную
соболевскую шкалу и позволяет значительно тоньше охарактеризовать регулярность
функций/распределений. Так, RO-меняющаяся функция может не иметь числового порядка
$s$ изменения на $\infty$.

В работах [14, 26~-- 31] найдены применения расширенной соболевской шкалы к
эллиптическим дифференциальным операторам, заданных в $\mathbb{R}^{n}$ и на
замкнутых гладких многообразиях. Отметим, что в последнее время пространства
Хермандера и их разные аналоги, именуемые пространствами обобщенной гладкости,
вызывают немалый интерес как сами по себе, так и с точки зрения приложений [10~-- 13].

Настоящая статья посвящена применениям расширенной соболевской шкалы к ЭКЗ, заданным
в ограниченной евклидовой области с гладкой границей. Наша цель~--- установить для
этой шкалы теоремы о разрешимости общей регулярной ЭКЗ и свойствах ее решений.

Работа состоит из шести пунктов. Пункт~1 --- введение. В п.~2 дано определение
регулярной ЭКЗ. Следующий п.~3 посвящен необходимым в работе пространствам
Хермандера, которые образуют расширенную соболевскую шкалу на $\mathbb{R}^{n}$,
евклидовой области и ее гладкой границе. Основные результаты статьи сформулированы в
п.~4. Это~--- теоремы о фредгольмовости операторов, соответствующих задаче в
расширенной соболевской шкале, порожденных ими изоморфизмах, априорной оценке
решений, их локальной регулярности. В качестве применений приведены достаточные
условия непрерывности обобщенных частных производных решений, а также классичности
обобщенного решения. Пункт~5 содержит необходимые для доказательства свойства
расширенной соболевской шкалы, в частности, интерполяционные. Все теоремы доказаны в
последнем п.~6.

\textbf{2. Постановка задачи.} Пусть $\Omega$ --- ограниченная область в евклидовом
пространстве $\mathbb{R}^{n}$, где  $n\geq2$. Предположим, что ее граница
$\Gamma:=\partial\Omega$ является бесконечно гладким замкнутым многообразием
(размерности $n-1$). Как обычно, $\overline{\Omega}:=\Omega\cup\Gamma$.

В области $\Omega$ рассматривается краевая задача
\begin{gather}\label{2f1}
Au(x)\equiv\sum_{|\mu|\leq{2q}}a_{\mu}(x)D^{\mu}u(x)=f(x),\quad x\in\Omega,\\
B_{j}u(x)\equiv \sum_{|\mu|\leq m_{j}}b_{j,\mu}(x)D^{\mu}u(x)=g_{j}(x),\quad
x\in\Gamma. \label{2f2}
\end{gather}
Здесь $A=A(x,D)$~--- линейное дифференциальное выражение на $\overline{\Omega}$
четного порядка $2q\geq2$, а каждое $B_{j}=B_{j}(x,D)$~--- граничное линейное
дифференциальные выражения на $\Gamma$ порядка $m_{j}\leq2q-1$. Все коэффициенты
дифференциальных выражений $A$ и $B_{j}$ предполагаются бесконечно гладкими
комплекснозначными функциями: $a_\mu \in C^{\infty}(\overline{\Omega})$ и
$b_{j,\mu}\in C^\infty(\Gamma)$. Положим $B:=(B_{1},\ldots,B_q)$.

Мы используем стандартные обозначения: $\mu:=(\mu_1,\ldots,\mu_n)$~--- мультииндекс,
$|\mu|:=\mu_1+\ldots+\mu_{n}$, $D^{\mu}:=D_1^{\mu_1}\ldots D_n^{\mu_n}$,
$D_k:=i\partial/\partial x_{k}$ при $k=1,\ldots,n$, где $i$~--- мнимая единица, а
$x=(x_1,\ldots,x_n)$~--- произвольная точка пространства $\mathbb{R}^{n}$.

Всюду в работе предполагается, что краевая задача \eqref{2f1}, \eqref{2f2} является
регулярной эллиптической  в области $\Omega$. Это означает что выражение $A$
правильно эллиптическое в $\overline{\Omega}$, а набор граничных выражений $B$
нормальный и удовлетворяет условию дополнительности по отношению к $A$ на $\Gamma$
(см., например, \cite{LionsMagenes71} (гл.~2, п.~1.4) или \cite{Triebel95}
(п.~5.2.1)). Из условия нормальности вытекает, что порядки $m_{j}$ граничных
дифференциальных выражений все различны.

Наряду с \eqref{2f1}, \eqref{2f2} рассмотрим краевую задачу
\begin{gather}\label{2f3}
A^{+}v\equiv\sum_{|\mu|\leq2q}D^{\mu}(\overline{a_{\mu}}\,v)=
\omega\quad\mbox{в}\quad\Omega,\\
B^{+}_{j}v=h_{j}\quad\mbox{на}\quad\Gamma,\quad j=1,\ldots,q. \label{2f4}
\end{gather}
Она формально сопряжена к задаче \eqref{2f1}, \eqref{2f2} относительно формулы Грина
\begin{equation*}
(Au,v)_{\Omega} + \sum^{q}_{j=1}(B_{j}u,C^{+}_{j}v)_{\Gamma} =
(u,A^{+}v)_{\Omega} + \sum_{j=1}^{q}(C_{j}u,B^{+}_{j}v)_{\Gamma},
\end{equation*}
справедливой для любых функций $u,v\in C^{\infty}(\overline{\Omega}).$ Здесь
$\{B^+_j\}$, $\{C_j\}$ и $\{C^+_j\}$~--- нормальные системы граничных линейных
дифференциальных выражений с коэффициентами класса $C^\infty(\Gamma)$. Их порядки
удовлетворяют условию
$$
\mathrm{ord}\,B_j+\mathrm{ord}\,C^+_j=\mathrm{ord}\,C_j+\mathrm{ord}\,B^+_j=2q-1.
$$
Здесь и далее через $(\cdot,\cdot)_{\Omega}$ и $(\cdot,\cdot)_{\Gamma}$ обозначены
скалярные произведения в комплексных гильбертовых пространствах $L_{2}(\Omega)$ и
$L_{2}(\Gamma)$, состоящих из функций, интегрируемых с квадратом в $\Omega$ и на
$\Gamma$ соответственно, а также расширения по непрерывности этих скалярных
произведений.

Краевая задача регулярная эллиптическая тогда и только тогда, когда формально
сопряженная к ней задача является регулярной эллиптической \cite{LionsMagenes71}
(гл.~2, п.~2.5).

Обозначим
\begin{gather*}
N:=\bigl\{u\in C^{\infty}(\,\overline{\Omega}\,):\;Au=0\;\;\mbox{в}\;\; \Omega,\;\;
B_{j}u=0\;\;\mbox{на}\;\;\Gamma\;\;\mbox{для}\;\;j=1,\ldots,q\bigr\},\\
N^{+}:=\bigl\{v\in C^{\infty}(\,\overline{\Omega}\,):\;A^{+}v=0\;\;\mbox{в}\;\;
\Omega,\;\;
B^{+}_{j}v=0\;\;\mbox{на}\;\;\Gamma\;\;\mbox{для}\;\;j=1,\ldots,q\bigr\}.
\end{gather*}
Поскольку задачи \eqref{2f1}, \eqref{2f2} и \eqref{2f3}, \eqref{2f4} являются
регулярными эллиптическими, пространства $N$ и $N^+$ конечномерны
\cite{LionsMagenes71} (гл.~2, п.~5.4). Отметим, что множество $N^{+}$ не зависит от
выбора сопряженной системы граничных выражений $\{B^{+}_{1},\ldots, B^{+}_{q}\}$,
удовлетворяющей формуле Грина \cite{LionsMagenes71} (гл.~2, п.~2.5).

Мы исследуем линейное отображение $u\mapsto(Au,Bu)$ в расширенной соболевской шкале.

\textbf{3. Расширенная  соболевская шкала} состоит из гильбертовых изотропных
пространств Хермандера $H^\varphi$, для которых показателем гладкости служит
произвольный функциональный параметр $\varphi\in\mathrm{RO}$. Приведем определение
класса $\mathrm{RO}$ и пространства~$H^\varphi$.

Класс $\mathrm{RO}$ состоит из всех измеримых по Борелю функций
$\varphi:[1,\infty)\rightarrow(0,\infty)$, для которых существуют числа $a>1$ и
$c\geq1$ такие, что
\begin{equation*}
c^{-1}\leq\frac{\varphi(\lambda t)}{\varphi(t)}\leq c\quad\mbox{для
любых}\quad t\geq1,\;\lambda\in[1,a]
\end{equation*}
(постоянные $a$ и $c$ могут зависеть от $\varphi\in\mathrm{RO}$). Такие функции
называют RO- (или OR-) меняющимися на бесконечности. Класс RO-меняющихся функций
введен В.~Г.~Авакумовичем \cite{Avakumovic36} и достаточно полно изучен (см.
\cite{Seneta76} (приложение~1) и \cite{BinghamGoldieTeugels89} (пп.~2.0~-- 2.2)).

Этот класс допускает простое описание:
$$
\varphi\in\mathrm{RO}\;\;\Leftrightarrow\;\;\varphi(t)=\exp\Biggl(\beta(t)+
\int\limits_{1}^{\:t}\frac{\gamma(\tau)}{\tau}\;d\tau\Biggr)\;\,\mbox{при}\;\,t\geq1,
$$
где действительны функции  $\beta$ и $\gamma$ измеримы по Борелю и ограничены на
полуоси $[1,\infty)$ (см., например, \cite{Seneta76} (приложение~1, теорема 1)).

Нам понадобится следующее свойство класса $\mathrm{RO}$ \cite{Seneta76}
(приложение~1, теорема 2). Для каждой функции $\varphi\in\mathrm{RO}$ существуют
числа $s_0,s_1\in\mathbb{R}$, $s_0\leq s_1$, и $c_0,c_1>0$ такие, что
\begin{equation}\label{2f5}
c_{0}\lambda^{s_{0}}\leq\frac{\varphi(\lambda t)}{\varphi (t)}\leq
c_{1}\lambda^{s_{1}} \quad\mbox{для всех}\quad t\geq1,\;\;\lambda\geq1.
\end{equation}
Положим
\begin{gather*}
\sigma_{0}(\varphi):=
\sup\,\{s_{0}\in\mathbb{R}:\,\mbox{верно левое неравенство в \eqref{2f5}}\},\\
\sigma_{1}(\varphi):=\inf\,\{s_{1}\in\mathbb{R}:\,\mbox{верно правое неравенство в
\eqref{2f5}}\};
\end{gather*}
здесь $-\infty<\sigma_{0}(\varphi)\leq\sigma_{1}(\varphi)<\infty$. Числа
$\sigma_{0}(\varphi)$ і $\sigma_{1}(\varphi)$ являются соответственно нижним и
верхним индексами Матушевской \cite{Matuszewska64} функции $\varphi\in\mathrm{RO}$
(см. \cite{BinghamGoldieTeugels89} (п.~2.1.2)).

Пусть $\varphi\in\mathrm{RO}$. По определению комплексное линейное пространство
$H^{\varphi}(\mathbb{R}^{n})$, где $n\geq1$, состоит из всех распределений
$w\in\mathcal{S}'(\mathbb{R}^{n})$ таких, что их преобразование Фурье
$\widehat{w}:=\mathcal{F}w$ локально суммируемо по Лебегу в $\mathbb{R}^{n}$ и
удовлетворяет условию
\begin{equation*}
\int\limits_{\mathbb{R}^{n}}\varphi^2(\langle\xi\rangle)\,|\widehat{w}(\xi)|^2\,d\xi
<\infty.
\end{equation*}
Здесь $\mathcal{S}'(\mathbb{R}^{n})$~--- линейное топологическое пространство
Л.~Шварца медленно растущих комплекснозначных распределений, заданных в
$\mathbb{R}^{n}$, а $\langle\xi\rangle:=(1+|\xi|^{2})^{1/2}$~--- сглаженный модуль
вектора $\xi\in\mathbb{R}^{n}$. Нам удобно трактовать распределения как
\emph{анти}линейные функционалы на пространстве $\mathcal{S}(\mathbb{R}^{n})$
основных функций.

В пространстве  $H^{\varphi}(\mathbb{R}^{n})$  определено скалярное
произведение распределений $w_1$, $w_2$ по формуле
\begin{equation*}
(w_1,w_2)_{H^{\varphi}(\mathbb{R}^{n})}:=
\int\limits_{\mathbb{R}^{n}}\varphi^2(\langle\xi\rangle)\,
\widehat{w_1}(\xi)\,\overline{\widehat{w_2}(\xi)}\,d\xi.
\end{equation*}
Оно задает на $H^{\varphi}(\mathbb{R}^{n})$ структуру гильбертового пространства и
определяет норму
\begin{equation*}
\|w\|_{H^{\varphi}(\mathbb{R}^{n})}:=(w,w)_{H^{\varphi}(\mathbb{R}^{n})}^{1/2}.
\end{equation*}
Это пространство сепарабельное; в нем плотно множество
$C_{0}^{\infty}(\mathbb{R}^{n})$ бесконечно дифференцируемых функций на
$\mathbb{R}^{n}$ с компактным носителем.

Пространство $H^{\varphi}(\mathbb{R}^{n})$~--- гильбертов изотропный случай
пространств $B_{p,k}$, введенных  и систематически изученных Л.~Хермандером в
\cite{Hermander63} (п.~2.2) (см. также его монографию \cite{Hermander83} (п.~10.1)).
А именно, $H^{\varphi}(\mathbb{R}^{n})=B_{p,k}$, если $p=2$ и
$k(\xi)=\varphi(\langle\xi\rangle)$ при $\xi\in\mathbb{R}^{n}$. Отметим, что при
$p=2$ пространства Хермандера совпадают с пространствами, введенными и изученными
Л.~Р.~Волевичем и Б.~П.~Панеяхом \cite{VolevichPaneah65} (\S~2).

Если функция $\varphi$ степенная: $\varphi(t)\equiv t^{s}$, то
$H^{\varphi}(\mathbb{R}^{n})=:H^{(s)}(\mathbb{R}^{n})$ является (гильбертовым)
пространством Соболева порядка $s\in\mathbb{R}$.

Вообще,
\begin{equation}\label{2f6}
s_{0}<\sigma_{0}(\varphi)\leq\sigma_{1}(\varphi)<s_{1}\;\;\Rightarrow\;\;
H^{(s_1)}(\mathbb{R}^{n})\hookrightarrow H^{\varphi}(\mathbb{R}^{n})\hookrightarrow
H^{(s_0)}(\mathbb{R}^{n}),
\end{equation}
причем оба вложения непрерывны и плотны.

Следуя \cite{13UMJ3}, класс гильбертовых функциональных пространств
\begin{equation}\label{2f7}
\{H^{\varphi}(\mathbb{R}^{n}):\varphi\in\mathrm{RO}\}
\end{equation}
называем расширенной соболевской шкалой на $\mathbb{R}^{n}$.

Нам необходимы ее аналоги для евклидовой области $\Omega$ и замкнутого многообразия
$\Gamma$. Эти аналоги строятся стандартным образом по классу \eqref{2f7} (см.
\cite{arXiv:1106.2049} (п.~2) и \cite{MikhailetsMurach10} (п.~2.4.2)). Приведем
соответствующие определения. Как и прежде, $\varphi\in\mathrm{RO}$.

Линейное пространство $H^{\varphi}(\Omega)$ состоит из сужений в область $\Omega$
всех распределений $w\in H^{\varphi}(\mathbb{R}^{n})$. В пространстве
$H^{\varphi}(\Omega)$ определена норма распределения  $u$ по формуле:
\begin{equation*}
\|u\|_{H^{\varphi}(\Omega)}:=\inf\bigl\{\,\|w\|_{H^{\varphi}(\mathbb{R}^{n})}:\,
w\in H^{\varphi}(\mathbb{R}^{n}),\;w=u\;\,\mbox{в}\;\,\Omega\,\bigr\}.
\end{equation*}
Относительно этой нормы пространство $H^{\varphi}(\Omega)$ гильбертово и
сеперабельно, поскольку оно является факторпространством гильбертова пространства
$H^\varphi(\mathbb{R}^{n})$ по подпространству
\begin{equation*}
\bigl\{w\in H^{\varphi}(\mathbb{R}^{n}):\,
\mathrm{supp}\,w\subseteq\mathbb{R}^{n}\setminus\Omega\bigr\}.
\end{equation*}
Множество  $C^\infty(\overline{\Omega})$ плотно в
$H^\varphi(\Omega)$.

Комплексное линейное пространство $H^{\varphi}(\Gamma)$ состоит из всех
распределений на $\Gamma$, которые в локальных координатах  принадлежат к
$H^{\varphi}(\mathbb{R}^{n-1})$. А именно, пусть произвольно выбран конечный атлас
из $C^{\infty}$-структуры на многообразии $\Gamma$, образованный локальными картами
$\alpha_j: \mathbb{R}^{n-1}\leftrightarrow U_{j}$, где $j=1,\ldots,\varkappa$.
(Открытые множества $U_{1},\ldots,U_{r}$ составляют покрытие многообразия $\Gamma$.)
Пусть, кроме того, произвольно выбраны функции $\chi_j\in C^{\infty}(\Gamma)$, где
$j=1,\ldots,\varkappa$, образующие разбиение единицы на $\Gamma$ и удовлетворяющие
условию $\mathrm{supp}\,\chi_j\subset U_j$. Тогда
\begin{equation*}
H^{\varphi}(\Gamma):=\bigl\{h\in\mathcal{D}'(\Gamma):\,
(\chi_{j}h)\circ\alpha_{j}\in H^{\varphi}(\mathbb{R}^{n-1})\;\;\mbox{для
каждого}\;\;j\in\{1,\ldots,\varkappa\}\bigr\}.
\end{equation*}
Здесь $\mathcal{D}'(\Gamma)$~--- линейное топологическое пространство всех
распределений на многообразии~$\Gamma$, а $(\chi_{j}h)\circ\alpha_{j}$~---
представление распределения $\chi_{j}u$ в локальной карте $\alpha_{j}$. В
пространстве  $H^\varphi (\Gamma)$ определено скалярное произведение распределений
$h_{1}$ и $h_{2}$ по формуле
$$
(h_{1},h_{2})_{H^{\varphi}(\Gamma)}:=
\sum_{j=1}^{\varkappa}\,((\chi_{j}h_{1})\circ\alpha_{j},
(\chi_{j}\,h_{2})\circ\alpha_{j})_{H^{\varphi}(\mathbb{R}^{n-1})}.
$$
Это пространство гильбертово сепарабельное и с точностью до эквивалентности норм не
зависит от выбора атласа и разбиения единицы \cite{MikhailetsMurach10}
(теорема~2.21). Множество $C^\infty(\Gamma)$ плотно в $H^\varphi(\Gamma)$.

Определеные выше функциональные пространства образуют расширенные соболевские шкалы
\begin{equation}\label{2f8}
\{H^{\varphi}(\Omega):\varphi\in\mathrm{RO}\}\quad\mbox{и}\quad
\{H^{\varphi}(\Gamma):\varphi\in\mathrm{RO}\}
\end{equation}
на $\Omega$ и $\Gamma$ соответственно. Они содержат шкалы гильбертовых пространств
Соболева: если $\varphi(t)\equiv t^s$ для некоторого $s\in\mathbb{R},$ то
$H^\varphi(\Omega)=:H^{(s)}(\Omega)$ и $H^\varphi(\Gamma)=:H^{(s)}(\Gamma)$ являются
пространствами Соболева порядка~$s$.

Необходимые нам свойства шкал \eqref{2f8} будут рассмотрены в п.~5.

\textbf{4. Основные результаты.} Сформулируем основные результаты статьи о свойствах
эллиптической краевой задачи \eqref{2f1}, \eqref{2f2} в расширенной соболевской
шкале. Их доказательства будут даны в п.~6.

Положим $\varrho(t):=t$ для произвольного $t\geq1$. Если $\varphi\in\mathrm{RO}$ и
$s\in\mathbb{R}$, то функция $\varrho^s\varphi\in\mathrm{RO}$, при этом ее индексы
Матушевской $\sigma_j(\varrho^s\varphi)=s+\sigma_j(\varphi)$ для каждого
$j\in\{0,1\}$.

\textbf{Теорема 1.} \it Пусть $\varphi\in\mathrm{RO}$ и $\sigma_0(\varphi)>-1/2$.
Тогда отображение $u\mapsto(Au,Bu)$, где $u\in C^\infty(\overline{\Omega})$,
продолжается единственным образом (по непрерывности) до ограниченного оператора
\begin{equation}\label{2f9}
(A,B):\,H^{\varphi\varrho^{2q}}(\Omega)\rightarrow
H^{\varphi}(\Omega)\oplus\bigoplus_{j=1}^{q}H^{\varphi\varrho^{2q-m_j-1/2}}(\Gamma)
=:\mathcal{H}^\varphi(\Omega,\Gamma).
\end{equation}
Этот оператор фредгольмов. Его ядро совпадает с $N$, а область значений состоит из
всех векторов $(f,g_1,\ldots,g_q)\in\mathcal{H}^\varphi(\Omega,\Gamma)$ таких, что
\begin{equation}\label{2f10}
(f,v)_{\Omega}+\sum_{j=1}^{q}\,(g_{j},C^{+}_{j}v)_{\Gamma}=0\quad \mbox{для
всех}\quad v\in N^{+}.
\end{equation}
Индекс оператора \eqref{2f9} равен $\dim N-\dim N^{+}$ и не зависит от $\varphi$.
\rm

Напомним, что линейный  оператор $T:X\rightarrow Y$, где $X$ и $Y$~--- банаховы
пространства, называется фредгольмовым, если его ядро $\ker T$ и коядро $Y/T(X)$
конечномерны. Если оператор $T$ фредгольмов, то его область значений $T(X)$ замкнута
в $Y$ (см., например, \cite{Hermander85} (лемма 19.1.1)), а индекс
$\mathrm{ind}\,T:=\dim\ker T-\dim(Y/T(X))$ конечен.

В частном случае, когда $N=\{0\}$ и $N^{+}=\{0\}$, оператор $(A,B)$ осуществляет
изоморфизм пространства $H^{\varphi\varrho^{2q}}(\Omega)$ на
$\mathcal{H}^\varphi(\Omega,\Gamma)$. В общем случае изоморфизм удобно строить с
помощью следующих проекторов.

Представим пространства, в которых действует оператор \eqref{2f9}, в виде прямых
сумм (замкнутых) подпространств:
\begin{gather}\label{2f11}
H^{\varphi\varrho^{2q}}(\Omega)=N\dotplus\bigl\{u\in
H^{\varphi\varrho^{2q}}(\Omega):\,(u,w)_\Omega=0\;\mbox{для всех}\;w\in N\bigr\},\\
\mathcal{H}^{\varphi}(\Omega,\Gamma)=\{(v,0,\ldots,0):\,v\in N^+\}
\dotplus(A,B)\bigl(H^{\varphi\varrho^{2q}}(\Omega)\bigr). \label{2f12}
\end{gather}
Такие разложения в прямые суммы существуют. В самом деле, \eqref{2f11} является
сужением на $H^{\varphi\varrho^{2q}}(\Omega)$ разложения пространства
$L_{2}(\Omega)$ в ортогональную сумму подпространств $N$ и $L_{2}(\Omega)\ominus N$.
Формула \eqref{2f12} верна, поскольку в ее правой части слагаемые имеют тривиальное
пересечение и (конечная) размерность первого слагаемого совпадает с коразмерностью
второго в силу теоремы~1.

Обозначим через $P$ и $P^+$ проекторы соответственно пространств
$H^{\varphi\varrho^{2q}}(\Omega)$ и $\mathcal{H}^{\varphi}(\Omega,\Gamma)$ на второе
слагаемое в суммах \eqref{2f11} и \eqref{2f12} параллельно первому слагаемому.
Отображения $P$ и $P^+$ не зависят от~$\varphi$.

 \textbf{Теорема 2.} \it Пусть $\varphi\in\mathrm{RO}$ и
$\sigma_0(\varphi)>-1/2$. Сужение отображения \eqref{2f9} на подпространство
$P(H^{\varphi\varrho^{2q}}(\Omega))$ является изоморфизмом
\begin{equation}\label{2f13}
(A,B):\,P(H^{\varphi\varrho^{2q}}(\Omega))\leftrightarrow P^+
(H^{\varphi}(\Omega,\Gamma)).
\end{equation} \rm

Для решения краевой задачи \eqref{2f1}, \eqref{2f2} выполняется следующая априорная
оценка.

\textbf{Теорема 3.} \it Пусть $\varphi\in\mathrm{RO}$ и $\sigma_0(\varphi)>-1/2$.
Тогда существует число $c=c(\varphi)>0$ такое, что
\begin{equation}\label{2f14}
\|u\|_{H^{\varphi\varrho^{2q}}(\Omega)}\leq c\,\bigl(\|(A,
B)u\|_{H^{\varphi}(\Omega, \Gamma)}+\|u\|_{L_2(\Omega)}\bigr)
\end{equation}
для произвольной функции $u\in H^{\varphi\varrho^{2q}}(\Omega).$ Здесь $c$ не
зависит от $u$. \rm

Исследуем локальную регулярность решения краевой задачи \eqref{2f1}, \eqref{2f2} в
расширенной соболевской шкале. Обозначим
$$
H^{2q-1/2+}(\Omega):=\bigcup_{l>2q-1/2}H^{(l)}(\Omega)=
\bigcup_{\substack{\alpha\in\mathrm{RO},\\\sigma_{0}(\alpha)>2q-1/2}}H^{\alpha}(\Omega)
$$
(последнее равенство верно ввиду \eqref{2f6}). Согласно теореме~1, для каждой
функции $u\in H^{2q-1/2+}(\Omega)$ определены по замыканию правые части краевой
задачи \eqref{2f1}, \eqref{2f2}.

Пусть $V$ --- произвольное открытое множество в $\mathbb{R}^{n}$ такое, что
$\Omega_0:=\Omega\cap V\neq\varnothing$. Положим $\Gamma_{0}:=\Gamma\cap V$
(возможный случай, когда $\Gamma_{0}=\varnothing$). Введем локальные аналоги
пространств $H^{\eta}(\Omega)$ и $H^{\eta}(\Gamma)$, где $\eta\in\mathrm{RO}$.
Обозначим
\begin{gather}\notag
H^{\alpha}_{\mathrm{loc}}(\Omega_{0},\Gamma_{0}):=\bigl\{u\in
\mathcal{D}'(\Omega):\,\chi u\in H^{\alpha}(\Omega)\;\\
\mbox{для всех}\;\chi\in C^{\infty}(\overline{\Omega})\;\mbox{таких, что}\;
\mathrm{supp}\,\chi\subset\Omega_0\cup\Gamma_{0}\bigr\}.\label{2f15}
\end{gather}
Здесь $\mathcal{D}'(\Omega)$~--- линейное топологическое пространство всех
распределений в области $\Omega$. Топология в линейном пространстве
$H^{\alpha}_{\mathrm{loc}}(\Omega_{0},\Gamma_{0})$ задается полунормами
$u\mapsto\|\chi u\|_{H^{\alpha}(\Omega)}$, где $\chi$~--- произвольная функция,
удовлетворяющая условию \eqref{2f15}. Аналогично определяем
$$
H^{\alpha}_{\mathrm{loc}}(\Gamma_{0}):=\bigl\{h\in D'(\Gamma):\,\chi h\in
H^{\alpha}(\Gamma)\;\mbox{для всех}\;\chi\in C^{\infty}(\Gamma)\;\mbox{таких, что}\;
\mathrm{supp}\,\chi\subset\Gamma_{0}\bigr\}.
$$
Топология в линейном пространстве $H^{\alpha}_{\mathrm{loc}}(\Gamma_{0})$ задается
полунормами $h\mapsto\|\chi h\|_{H^{\alpha}(\Gamma)}$, где $\chi\in
C^{\infty}(\Gamma)$ и $\mathrm{supp}\,\chi\subset\Gamma_{0}$.

\textbf{Теорема 4.} \it Предположим, что функция $u\in H^{2q-1/2+}(\Omega)$ является
решением краевой задачи \eqref{2f1}, \eqref{2f2}, правые части которой удовлетворяют
условиям
$$
f\in H^{\varphi}_{\mathrm{loc}}(\Omega_{0},\Gamma_{0})\quad\mbox{и}\quad g_{j}\in
H^{\varphi\varrho^{2q-m_{j}-1/2}}_{\mathrm{loc}}(\Gamma_{0})\;\;
\mbox{при}\;\;j=1,\ldots,q
$$
для некоторого функционального параметра $\varphi\in\mathrm{RO}$ такого, что
$\sigma_0(\varphi)>-1/2$. Тогда $u\in
H^{\varphi\varrho^{2q}}_{\mathrm{loc}}(\Omega_0,\Gamma_0)$. \rm

Отметим важные частные случаи этой теоремы.

Если $\Omega_0=\Omega$ и $\Gamma_0=\Gamma$, то
$H^{\alpha}_{\mathrm{loc}}(\Omega_{0},\Gamma_{0})=H^{\alpha}(\Omega)$ и
$H^{\alpha}_{\mathrm{loc}}(\Gamma_{0})=H^{\alpha}(\Gamma)$ для любого
$\alpha\in\mathrm{RO}$. В этом случае теорема~4 утверждает, что регулярность решения
$u$ повышается глобально, т.~е. во всей области $\Omega$ вплоть до ее
границы~$\Gamma$.

Если $\Omega_0=\Omega$ и $\Gamma_0=\varnothing$, то по теореме 4 регулярность
решения $u$ повышается в окрестности каждой (внутренней) точки области
$\overline{\Omega}$.

Теоремы 1~-- 4 распространяют на расширенную соболевскую шкалу известные свойства регулярных эллиптических краевых задач в пространствах Соболева \cite{Berezansky68, LionsMagenes71, Roitberg96, Agranovich97}, а также в уточненной соболевской шкале \cite{MikhailetsMurach10}.

В связи с теоремой~2 отметим, что в статье Г.~Шлензак \cite{Slenzak74}, установлен аналог изоморфизма \eqref{2f13} для некоторого класса гильбертовых пространств Хермандера. Однако, использованный Г.~Шлензак класс не имеет конструктивного описания.

В качестве применения расширенной соболевской шкалы приведем достаточное условие
непрерывности обобщенных частных производных решения эллиптической краевой задачи.

\textbf{Теорема 5.} \it Пусть выполняется условие теоремы~$4$. Кроме того,
предположим, что
\begin{equation}\label{2f16}
\int\limits_1^{\infty}t^{2k+n-1-4q}\varphi^{-2}(t)dt<\infty
\end{equation}
для некоторого целого $k\geq0$. Тогда $u\in C^{k}(\Omega_{0}\cup\Gamma_{0})$. \rm

\textbf{\textit{Замечание}~1.} В теореме~5 условие \eqref{2f16} не только
достаточное, но и необходимое на классе всех решений $u$, удовлетворяющих условию
теоремы~4.

Особо выделим достаточное условие классичности решения $u\in H^{2q-1/2+}(\Omega)$
краевой задачи \eqref{2f1}, \eqref{2f2}.  Положим $m:=\max\{m_{1},\ldots,m_{q}\}$.

\textbf{Теорема 6.} \it Предположим, что функция $u\in H^{2q-1/2+}(\Omega)$ является
решением краевой задачи \eqref{2f1}, \eqref{2f2}, правые части которой удовлетворяют
условиям
\begin{gather}\label{2f17}
f\in H^{\varphi_1}_{\mathrm{loc}}(\Omega,\varnothing)\cap H^{\varphi_2}(\Omega),\\
g_j\in H^{\varphi_2\varrho^{2q-m_j-1/2}}(\Gamma)\quad\mbox{при}\quad j=1,\ldots,q
\label{2f18}
\end{gather}
для некоторых функциональных параметров $\varphi_1,\varphi_2\in\mathrm{RO}$ таких,
что $\sigma_0(\varphi_1),\sigma_0(\varphi_2)>-1/2$, и
\begin{gather}\label{2f19}
\int\limits_1^{\infty}t^{n-1}\varphi_{1}^{-2}(t)dt<\infty,\\
\int\limits_1^{\infty}t^{2m+n-1-4q}\varphi_{2}^{-2}(t)dt<\infty. \label{2f20}
\end{gather}
Тогда $u\in C^{2q}(\Omega)\cap C^{m}(\overline{\Omega})$, т.~е. решение $u$
классическое. \rm

В связи с последней теоремой отметим следующее. Если $u$~--- классическое решение
рассматриваемой краевой задачи, то левые части равенств \eqref{2f1} и \eqref{2f2}
вычисляется с помощью классических производных, а сами равенства выполняются в
каждой точке множества $\Omega$ либо $\Gamma$ соответственно. При этом
\begin{equation}\label{2f21}
F=(A,B)u\in C(\Omega)\times(C(\Gamma))^{q}
\end{equation}
Заметим, что обратное утверждение, вообще говоря, не верно. А именно, из условия
\eqref{2f21} не следует, что решение $u$ является классическим (см., например,
\cite{GilbargTrudinger98} (п.~4.5, примечания)).

\textbf{5. Вспомогательные результаты.} Здесь мы изложим несколько полезных
результатов, которые будут использованы в доказательствах теорем. Первые два из них
относятся к свойствам вложений пространств Хермандера.

\textbf{Предложение 1.} \it Пусть $\alpha,\alpha_1\in\mathrm{RO}$. Отношение функций
$\alpha/\alpha_1$ ограничено в  окрестности $\infty$ тогда и только тогда, когда
$H^{\alpha_1}(\Omega)\hookrightarrow H^\alpha(\Omega)$. Это вложение непрерывно и
плотно. Оно компактно тогда и только тогда, когда
$\alpha(t)/\alpha_1(t)\rightarrow0$ при $t\rightarrow\infty$. Это предложение
сохраняет силу, если в нем заменить $\Omega$ на $\Gamma$. \rm

Предложение~1 является прямым следствием теорем 2.2.2 и 2.2.3 из монографии
Л.~Хермандера \cite{Hermander63}.

В силу этого предложения верны следующие вложения для пространств Хермандера и
Соболева: если $\alpha\in\mathrm{RO}$, $s_0<\sigma_0(\alpha)$ и $s_1
>\sigma_1(\alpha)$, то
\begin{gather}\label{2f22}
H^{(s_1)}(\Omega)\hookrightarrow H^\alpha(\Omega)\hookrightarrow
H^{(s_0)}(\Omega),\\
H^{(s_1)}(\Gamma)\hookrightarrow H^\alpha(\Gamma)\hookrightarrow H^{(s_0)}(\Gamma).
\label{2f23}
\end{gather}
Эти вложения компактны и плотны.

\textbf{Предложение 2.} \it Пусть произвольно выбраны функция $\alpha\in\mathrm{RO}$
и целое число $k\geq0$. Тогда
$$
\int\limits_1^{\infty} t^{2k+n-1}\alpha^{-2}(t)\,dt<\infty\;\;\Leftrightarrow\;\;
H^\alpha(\Omega)\hookrightarrow C^k(\overline{\Omega}).
$$
Последнее вложение непрерывно. \rm

Это предложение следует из теоремы 2.2.7 монографии Л.~Хермандера~\cite{Hermander63}
(см. соответствующие рассуждения в \cite{12UMJ11} (лемма~2)).

Связь между расширенной соболевской шкалой и пространствами Соболева не
исчерпывается вложениями \eqref{2f6} и их  аналогами \eqref{2f22} и \eqref{2f23};
она значительно глубже. А именно, каждое пространство $H^{\alpha}$, где
$\alpha\in\mathrm{RO}$, есть результат интерполяции с подходящим функциональным
параметром пары соболевских пространств $H^{(s_0)}$ и $H^{(s_1)}$, если
$s_{0}<\sigma_{0}(\alpha)$ и $\sigma_{1}(\alpha)<s_{1}$. Это фундаментальное
свойство будет использовано нами в доказательстве теоремы~1.

Метод интерполяции с функциональным параметром нормированных пространств впервые появился в статье
К.~Фояша, Ж.-Л.~Лионса \cite{FoiasLions61} и был исследован рядом авторов (см. монографии \cite{Ovchinnikov84, BrudnyiKrugljak91} и цитированную там литературу).

Для наших целей достаточно воспользоваться интерполяцией с функциональным параметром сепарабельных гильбертовых пространств. Напомним ее определение и некоторые свойства (см. монографию
\cite{MikhailetsMurach10} (пп. 1.1, 2.4.2) или статью \cite{08MFAT1} (п.~2), где эти
вопросы изложены систематически).

Пусть задана упорядоченная пара $X:=[X_{0},X_{1}]$ сепарабельных комплексных
гильбертовых пространств $X_{0}$ и $X_{1}$ такая, что выполняется непрерывное и
плотное вложение $X_{1}\hookrightarrow X_{0}$. Пару $X$ называем допустимой. Для нее
существует изометрический изоморфизм $J:X_{1}\leftrightarrow X_{\,0}$ такой, что $J$
является самосопряженным положительно определенным оператором в пространстве $X_{0}$
с областью определения $X_{1}$. Оператор $J$ определяется парой  $X$ однозначно; он
называется порождающим для $X$.

Обозначим  через $\mathcal{B}$ множество всех измеримых по Борелю функций
$\psi:\nobreak(0,\infty)\rightarrow(0,\infty)$, которые ограничены на каждом отрезке
$[a,b]$ и отделены от нуля на каждом множестве $[r,\infty)$, где $0<a<b<\infty$ и
$r>0$.

Пусть $\psi\in\mathcal{B}$. В пространстве $X_{0}$ определен, как функция от $J$,
оператор $\psi(J)$. Обозначим через $[X_{0},X_{1}]_\psi$ или, короче, $X_{\psi}$
область определения оператора $\psi(J)$, наделенную скалярным произведением
$$
(w_1, w_2)_{X_\psi}:=(\psi(J)w_1,\psi(J)w_2)_{X_0}
$$
и соответствующей нормой $\|w\|_{X_\psi}=(w,w)_{X_\psi}^{1/2}$. Пространство
$X_\psi$ гильбертово и сепарабельно, при этом выполняется непрерывное и плотное
вложение $X_\psi \hookrightarrow X_0$.

Функцию $\psi\in\mathcal{B}$ называем интерполяционным параметром,
если для произвольных допустимых пар $X=[X_0, X_1]$ и $Y=[Y_0, Y_1]$
гильбертовых пространств и для любого линейного отображения  $T$,
заданного на $X_0$, выполняется следующее. Если при каждом
$j\in\{0,1\}$ сужение отображения $T$ на пространство $X_{j}$
является ограниченным оператором $T:X_{j}\rightarrow Y_{j}$, то и
сужение отображения $T$ на пространство $X_\psi$ является
ограниченным оператором  $T:X_{\psi}\rightarrow Y_{\psi}$. Тогда
будем говорить, что пространство $X_\psi$ получено интерполяцией с
функциональным параметром $\psi$ пары $X$.

Функция $\psi\in\mathcal{B}$ является интерполяционным параметром тогда и только
тогда, когда она псевдовогнута  в окрестности бесконечности, то есть
$\psi(t)\asymp\psi_{1}(t)$ при $t\gg1$ для некоторой положительной  вогнутой функции
$\psi_{1}(t)$. (Здесь $\psi\asymp\psi_{1}$ означает ограниченность обоих отношений
$\psi/\psi_{1}$ и $\psi_{1}/\psi$ на указанном множестве). Этот важный факт вытекает
из теоремы Ж.~Петре \cite{Peetre68} об описании всех интерполяционных функций
положительной степени (см. также монографию \cite{BerghLefstrem76} (п.~5.4)).

Упомянутое выше интерполяционное свойство расширенной соболевской
шкалы формулируется следующим образом

\textbf{Предложение 3.} \it Пусть заданы функция $\alpha\in\mathrm{RO}$ и числа
$s_0,s_1\in\mathbb{R}$ такие, что $s_0<\sigma_0(\alpha)$ и $s_1>\sigma_1(\alpha)$.
Положим
\begin{equation}\label{2f24}
\psi(t)=
\begin{cases}
\;t^{{-s_0}/{(s_1-s_0)}}\,
\alpha\bigl(t^{1/{(s_1-s_0)}}\bigr)&\text{при}\quad t\geq1, \\
\;\alpha(1)&\text{при}\quad0<t<1.
\end{cases}
\end{equation}
Тогда функция $\psi\in\mathcal{B}$ является интерполяционным параметром, и
выполняются следующие равенства пространств с эквивалентностью  норм в них:
\begin{gather}\label{2f25}
\bigl[H^{(s_0)}(\Omega),H^{(s_1)}(\Omega)\bigr]_{\psi}=H^{\alpha}(\Omega),\\
\bigl[H^{(s_0)}(\Gamma),H^{(s_1)}(\Gamma)\bigr]_{\psi}=H^{\alpha}(\Gamma).
\label{2f26}
\end{gather} \rm

Формула \eqref{2f25} установлена в \cite{arXiv:1106.2049} (теорема 5.1), а
\eqref{2f26}~--- в \cite{MikhailetsMurach10} (теорема 2.21). Отметим, что при
условии предложения~3 также
$$
\bigl[H^{(s_0)}(\mathbb{R}^{n}),H^{(s_1)}(\mathbb{R}^{n})\bigr]_{\psi}=
H^{\alpha}(\mathbb{R}^{n})
$$
с равенством норм \cite{MikhailetsMurach10} (теорема 2.19).

Укажем еще два важных свойства расширенной соболевской шкалы
$\{H^{\alpha}:\alpha\in\mathrm{RO}\}$. Она замкнута относительно интерполяции с
функциональным  параметром и совпадает (с точностью до эквивалентности норм) с
классом всех гильбертовых пространств, интерполяционных  для пар соболевских
пространств $[H^{(s_0)},H^{(s_1)}]$, где $s_0,s_1\in\mathbb{R}$ и $s_0<s_1$.
Последнее свойство следует из теоремы В.~И.~Овчинникова \cite{Ovchinnikov84}
(п.~11.4) об описании всех гильбертовых пространств, интерполяционных для заданной
пары гильбертовых пространств (см. пояснения в \cite{MikhailetsMurach10} (п.
2.4.2)). Напомним, что свойство (гильбертового) пространства $H$ быть
интерполяционным для допустимой пары $X=[X_0,X_1]$ означает следующее:
а)~выполняются непрерывные вложения $X_1\hookrightarrow H\hookrightarrow X_0$,
б)~любой линейный оператор, ограниченный на каждом из пространств $X_0$ и $X_1$,
является ограниченным и на~$X$.

В конце этого пункта приведем два общих свойства интерполяции, которые будут
использованы  в доказательствах. Первое из них показывает, что при интерполяции
пространств наследуется не только ограниченность, но и фредгольмовость линейных
операторов при некоторых дополнительных условиях \cite{MikhailetsMurach10}
(п.~1.1.7).

\textbf{Предложение 4.} \it Пусть $X=[X_0,X_1]$ и $Y=[Y_0,Y_1]$ --- допустимые пары
гильбертовых пространств. Пусть, кроме того, на $X_0$ задано линейное отображение
$T$ такое, что его сужение на пространства $X_j$, где $j=0,1$, являются
фредгольмовыми ограниченными операторами $T:X_j\rightarrow Y_j$, которые имеют общее
ядро и одинаковый индекс. Тогда для произвольного интерполяционного параметра
$\psi\in\mathcal{B}$ ограниченный оператор $T:X_\psi\rightarrow Y_\psi$ фредгольмов
с тем же  ядром и индексом, а его  область значений равна $Y_\psi\cap T(X_0)$. \rm

Второе свойство сводит интерполяцию ортогональных сумм гильбертовых пространств к
интерполяции слагаемых \cite{MikhailetsMurach10} (п.~1.1.5).

\textbf{Предложение 5.} \it Пусть задано конечное число допустимых пар
$[X_{0}^{(k)},X_{1}^{(k)}]$  гильбертовых пространств, где $k=1,\ldots,p$. Тогда для
произвольного $\psi\in\mathcal{B}$ выполняется
\begin{equation*}
\biggl[\,\bigoplus_{k=1}^{p}X_{0}^{(k)},\,\bigoplus_{k=1}^{p}X_{1}^{(k)}\biggr]_{\psi}=\,
\bigoplus_{k=1}^{p}\bigl[X_{0}^{(k)},\,X_{1}^{(k)}\bigr]_{\psi}
\end{equation*}
с равенством  норм. \rm

\textbf{6. Доказательство основных результатов.}

\textbf{\textit{Доказательство теоремы}~1\textit{.}} В случае пространств Соболева,
т.~е. когда $\varphi(t)\equiv\nobreak t^{s}$ и  $s>-1/2$, эта теорема известна (см.
\cite{LionsMagenes71} (гл.~2, теорема~5.4) для $s\geq0$ и \cite{Roitberg96} (теоремы
4.1.3 и 5.3.1) для $s>-1/2$). В общей ситуации произвольного $\varphi\in\mathrm{RO}$
мы выведем теорему~1 из соболевского случая с помощью интерполяции с функциональным
параметром.

Выберем действительные числа $l_{0}$ и $l_{1}$ так, чтобы
$-1/2<l_{0}<\sigma_{0}(\varphi)$ и $\sigma_{1}(\varphi)<l_{1}$. Отображение
$u\mapsto(Au,Bu)$, где $u\in C^\infty(\overline{\Omega})$, продолжается по
непрерывности до фредгольмовых ограниченных операторов
\begin{equation}\label{2f27}
(A,B):\,H^{(l_{r}+2q)}(\Omega)\rightarrow
H^{(l_{r})}(\Omega)\oplus\bigoplus_{j=1}^{q}
H^{(l_{r}+2q-m_{j}-1/2)}(\Gamma)=:\mathcal{H}^{(l_{r})}(\Omega,\Gamma), \quad
r=0,\,1,
\end{equation}
действующих в пространствах Соболева. Эти операторы имеют общее ядро $N$ и
одинаковый индекс, равный $\dim N-\dim N^+$. Кроме того,
\begin{equation}\label{2f28}
(A,B)\bigl(H^{(l_{r}+2q)}(\Omega)\bigr)=
\bigl\{(f,g_1,\ldots,g_q)\in\mathcal{H}^{(l_{r})}(\Omega,\Gamma):\,\mbox{верно
\eqref{2f10}}\bigr\}.
\end{equation}

Определим функцию $\psi$ по формуле \eqref{2f24}, в которой полагаем
$\alpha:=\varphi$, $s_{0}:=l_{0}$ и $s_{1}:=l_{1}$. Согласно предложению~3 функция
$\psi$ является интерполяционным параметром. Поэтому в силу предложения~4 из
ограниченности и фредгольмовости обоих операторов \eqref{2f27} следует
ограниченность и фредгольмовость оператора
\begin{equation}\label{2f29}
(A,B):\,\bigl[H^{(l_{0}+2q)}(\Omega),H^{(l_{1}+2q)}(\Omega)\bigr]_{\psi}\rightarrow
\bigl[\mathcal{H}^{(l_{0})}(\Omega,\Gamma),
\mathcal{H}^{(l_{1})}(\Omega,\Gamma)\bigr]_{\psi}.
\end{equation}
Он является сужением  оператора \eqref{2f27}, где $r=0$. Покажем, что отображение
\eqref{2f29} есть оператор \eqref{2f9} из формулировки теоремы~1.

На основании предложения~3 имеем следующие равенства пространств с точностью до
эквивалентности  норм в них:
\begin{gather}\label{2f30}
\bigl[H^{(l_{0})}(\Omega),H^{(l_{1})}(\Omega)\bigr]_{\psi}=H^{\varphi}(\Omega),\\
\bigl[H^{(l_{0}+2q)}(\Omega),H^{(l_{1}+2q)}(\Omega)\bigr]_{\psi}=
H^{\varphi\varrho^{2q}}(\Omega);\label{2f31}\\
\bigl[H^{(l_{0}+2q-m_{j}-1/2)}(\Gamma),
H^{(l_{1}+2q-m_{j}-1/2)}(\Gamma)\bigr]_{\psi}=
H^{\varphi\varrho^{2q-m_{j}-1/2}}(\Gamma)\quad\mbox{при}\quad
j=1,\ldots,q.\label{2f32}
\end{gather}
Относительно двух последних формул отметим следующее. Равенство \eqref{2f31}
получаем, полагая  $\alpha:=\varphi\varrho^{2q}$, $s_{0}:=l_{0}+2q$ и
$s_{1}:=l_{1}+2q$ в предложении~3, а равенство \eqref{2f32} получаем, полагая там
$$
\alpha:=\varphi\varrho^{2q-m_{j}-1/2},\quad s_{0}:=l_{0}+2q-m_{j}-1/2,\quad
s_{1}:=l_{1}+2q-m_{j}-1/2.
$$
При этом для обеих формул функция $\psi$ удовлетворяет \eqref{2f24}.

Из формул \eqref{2f30} и \eqref{2f32} следует в силу предложения~5, что
\begin{equation}\label{2f33}
\begin{gathered}
\bigl[\mathcal{H}^{(l_{0})}(\Omega,\Gamma),
\mathcal{H}^{(l_{1})}(\Omega,\Gamma)\bigr]_{\psi}=
\bigl[H^{(l_{0})}(\Omega),H^{(l_{1})}(\Omega)\bigr]_{\psi}\oplus\\
\oplus\bigoplus_{j=1}^{q}\bigl[H^{(l_{0}+2q-m_{j}-1/2)}(\Gamma),
H^{(l_{1}+2q-m_{j}-1/2)}(\Gamma)\bigr]_{\psi}=\mathcal{H}^\varphi(\Omega,\Gamma)
\end{gathered}
\end{equation}
с точностью до эквивалентности норм.

Теперь на основании равенств \eqref{2f31} и \eqref{2f33} делаем вывод, что
ограниченный и фредгольмов оператор \eqref{2f29} действует в паре пространств
$(A,B):\,H^{\varphi\varrho^{2q}}(\Omega)\rightarrow\mathcal{H}^\varphi(\Omega,\Gamma)$.
В~силу предложения~4 ядро этого оператора и его индекс совпадают с общим ядром $N$ и
одинаковым индексом $\dim N-\dim N^+$ операторов \eqref{2f27}. Кроме того,
\begin{gather*}
(A,B)\bigl(H^{\varphi\varrho^{2q}}(\Omega)\bigr)=
\mathcal{H}^{\varphi}(\Omega,\Gamma)\cap(A,B)\bigl(H^{(l_{0}+2q)}(\Omega)\bigr)=\\
=\bigl\{(f,g_1,\ldots,g_q)\in\mathcal{H}^{\varphi}(\Omega,\Gamma):\,\mbox{верно
\eqref{2f10}}\bigr\}.
\end{gather*}
Здесь также воспользовались равенством \eqref{2f28} и вложением
$\mathcal{H}^\varphi(\Omega,\Gamma)\hookrightarrow\mathcal{H}^{(l_{0})}(\Omega,\Gamma)$,
которое следует из формул \eqref{2f22} и \eqref{2f23}. Остается отметить, что этот
оператор является продолжением по непрерывности отображения $u\mapsto(Au,Bu)$, где
$u\in C^\infty(\overline{\Omega})$, поскольку множество
$C^\infty(\overline{\Omega})$ плотно в пространстве
$H^{\varphi\varrho^{2q}}(\Omega)$.

Теорема 1 доказана.

\textbf{\textit{Доказательство теоремы}~2\textit{.}} По теореме~1 сужение оператора
\eqref{2f9} на $P(H^{\varphi\varrho^{2q}}(\Omega))$ является непрерывным и
биективным отображением подпространства $P(H^{\varphi\varrho^{2q}}(\Omega))$ на
подпространство $P^+(\mathcal{H}^{\varphi}(\Omega,\Gamma))$. Следовательно, по
теореме Банаха о обратном операторе, оно является изоморфизмом \eqref{2f13}.
Теорема~2 доказана.

\textbf{\textit{Доказательство теоремы}~3\textit{.}} Для произвольной функции $u\in
H^{\varphi\varrho^{2q}}(\Omega)$ в силу теоремы~2 имеем
\begin{gather*}
\|u\|_{H^{\varphi\varrho^{2q}}(\Omega)}\leq
\|Pu\|_{H^{\varphi\varrho^{2q}}(\Omega)}+\|u-Pu\|_{H^{\varphi\varrho^{2q}}(\Omega)}\leq\\
\leq
c_{1}\|(A,B)Pu\|_{\mathcal{H}^{\varphi}(\Omega,\Gamma)}+c_{2}\|u-Pu\|_{L_2(\Omega)}.
\end{gather*}
Здесь $c_1$~---  норма оператора, обратного к изоморфизму \eqref{2f13}, а $c_2$~---
некоторое положительное  число, не зависящее от $u$. Это число существует, поскольку
функция $u-Pu$ принадлежит конечномерному пространству $N$, а в нем эквивалентны все
нормы, в частности, нормы в пространствах $H^{\varphi\varrho^{2q}}(\Omega)$ и
$L_2(\Omega)$. Отсюда с учетом формул
$$
(A,B)Pu=(A,B)u\quad\mbox{і}\quad\|u-Pu\|_{L_2(\Omega)}\leq\|u\|_{L_2(\Omega)}
$$
имеем требуемую оценку \eqref{2f14}. Теорема~3 доказана.

\textbf{\textit{Доказательство теоремы}~4\textit{.}} По условию, $u\in
H^{(s+2q)}(\Omega)$ для некоторого числа $s$, удовлетворяющего неравенству
$-1/2<s<\sigma_{0}(\varphi)$. Отметим, что в силу \eqref{2f22} пространство
$H^{(s+2q)}(\Omega)$ шире, чем $H^{\varphi\varrho^{2q}}(\Omega)$.

Сначала установим теорему в случае глобальной регулярности, когда
$\Omega_{0}=\Omega$ и $\Gamma_{0}=\Gamma$. Тогда, по условию,
$F:=(A,B)u\in\mathcal{H}^{\varphi}(\Omega,\Gamma)$. Следовательно, на основании
теоремы~1 (использованной также в соболевском случае $\varphi(t)\equiv t^{s}$)
запишем
$$
F\in\mathcal{H}^{\varphi}(\Omega,\Gamma)\cap(A,B)\bigl(H^{(s+2q)}(\Omega)\bigr)=
(A,B)\bigl(H^{\varphi\varrho^{2q}}(\Omega)\bigr).
$$
Значит, наряду с условием $(A,B)u=F$ выполняется равенство $(A,B)v=F$ для некоторой
функции $v\in H^{\varphi\varrho^{2q}}(\Omega)$. Тогда $(A,B)(u-v)=0$, что по
теореме~1 (для $\varphi(t)\equiv t^{s}$) влечет за собой включения
$$
w:=u-v\in N\subset C^{\infty}(\overline{\Omega})\quad\mbox{и}\quad u=v+w\in
H^{\varphi\varrho^{2q}}(\Omega).
$$
В рассмотренном случае теорема~4 доказана.

Отсюда выведем ее в общем случае. Обозначим
$$
\Upsilon:=\bigl\{\chi\in
C^{\infty}(\overline{\Omega}):\,\mathrm{supp}\,\chi\subset\Omega_0\cup\Gamma_0\bigr\}.
$$
Сначала докажем, что в силу условия теоремы~4, выполняется импликация для каждого
$i\in\mathbb{N}$:
\begin{equation}\label{2f34}
\begin{gathered}
\bigl(\chi u\in H^{\varphi\varrho^{2q}}(\Omega)+
H^{(s+i-1+2q)}(\Omega)\;\;\mbox{для всех}\;\;\chi\in\Upsilon\bigr)\;\Rightarrow\\
\Rightarrow\;\bigl(\chi u\in H^{\varphi\varrho^{2q}}(\Omega)+
H^{(s+i+2q)}(\Omega)\;\;\mbox{для всех}\;\;\chi\in\Upsilon\bigr).
\end{gathered}
\end{equation}
Здесь и далее в доказательстве используем алгебраические суммы пространств.

Произвольно выберем $i\in\mathbb{N}$  и предположим, что истинна посылка импликации
\eqref{2f34}. Рассмотрим любую функцию $\chi\in\Upsilon$ и выберем функцию
$\eta\in\Upsilon$ такую, что $\eta=1$ в окрестности $\mathrm{supp}\,\chi$. По
условию, $\chi F\in\mathcal{H}^{\varphi}(\Omega,\Gamma)$, где $F:=(A,B)u$.

Переставив  оператор умножения  на функцию $\chi$  с дифференциальным оператором
$(A,B)$, можем записать следующее:
\begin{gather}\notag
\chi F=\chi(A,B)(\eta u)=(A,B)(\chi\eta u)-(A',B')(\eta u),\\
(A,B)(\chi u)=\chi F+(A',B')(\eta u). \label{2f35}
\end{gather}
Здесь $(A',B')$~--- некоторый дифференциальный оператор вида $(A,B)$, порядки
компонент которого меньше (по крайней мере на единицу), чем порядки соответствующих
компонент оператора $(A,B)$. В силу посылки импликации \eqref{2f34} имеем $\eta
u=u_{1}+u_{2}$ для некоторых функций
$$
u_{1}\in H^{\varphi\varrho^{2q}}(\Omega)\quad\mbox{и}\quad u_{2}\in
H^{(s+i-1+2q)}(\Omega).
$$
Отсюда и в силу \eqref{2f35} можем записать $(A,B)(\chi u)=F_{1}+F_{2}$, где
\begin{gather}\label{2f36}
F_{1}:=\chi F+(A',B')u_{1}\in\mathcal{H}^{\varphi}(\Omega,\Gamma),\\
F_{2}:=(A',B')u_{2}\in\mathcal{H}^{(s+i)}(\Omega,\Gamma). \label{2f37}
\end{gather}

Объясним включения, фигурирующие в формулах \eqref{2f36} и \eqref{2f37}. Здесь
$\mathcal{H}^{(l)}(\Omega,\Gamma)$ обозначает пространство
$\mathcal{H}^{\alpha}(\Omega,\Gamma)$ в соболевском случае $\alpha(t)\equiv t^{l}$ и
$l\in\mathbb{R}$. Поскольку  компонентами
$\mathrm{ord}\,(A',B')\leq\mathrm{ord}\,(A,B)-1$, то отображение
$v\mapsto(A'v,B'v)$, где $v\in\nobreak C^\infty(\overline{\Omega})$, продолжается по
непрерывности до ограниченного оператора
$$
(A',B'):\,H^{(l+2q)}(\Omega)\rightarrow\mathcal{H}^{(l+1)}(\Omega,\Gamma)\quad
\mbox{для каждого}\quad l>m-2q+1/2.
$$
Отсюда при $l:=s+i-1$ и из включения $u_{2}\in H^{(s+i-1+2q)}(\Omega)$ следует
\eqref{2f37}.

Далее, из ограниченности операторов
$$
(A',B'):\,H^{(l_{i}+2q)}(\Omega)\rightarrow
\mathcal{H}^{(l_{i}+1)}(\Omega,\Gamma)\hookrightarrow
\mathcal{H}^{(l_{i})}(\Omega,\Gamma)\quad\mbox{при}\quad i=0,\,1,
$$
следует в силу интерполяционных формул \eqref{2f31} и \eqref{2f33} ограниченность
оператора
\begin{equation}\label{2f38}
\begin{gathered}
(A',B'):\,H^{\varphi\varrho^{2q}}(\Omega)=
\bigl[H^{(l_{0}+2q)}(\Omega),H^{(l_{1}+2q)}(\Omega)\bigr]_{\psi}\rightarrow\\
\rightarrow\bigl[\mathcal{H}^{(l_{0})}(\Omega,\Gamma),
\mathcal{H}^{(l_{1})}(\Omega,\Gamma)\bigr]_{\psi}=
\mathcal{H}^{\varphi}(\Omega,\Gamma).
\end{gathered}
\end{equation}
Здесь числа $l_{0}$ и $l_{1}$ и интерполяционный параметр $\psi$ такие, как в
доказательстве теоремы~1. Теперь формула \eqref{2f36} является следствием
\eqref{2f38} и включений $\chi F\in\mathcal{H}^{\varphi}(\Omega,\Gamma)$, $u_{1}\in
H^{\varphi\varrho^{2q}}(\Omega)$.

Воспользуемся проектором $P^+$ и теоремой~2 (в соболевском случае тоже). Из
равенства $(A,B)(\chi u)=F_{1}+F_{2}$ и включений \eqref{2f36} и \eqref{2f37}
следует, что
$$
(A,B)(\chi u)=P^+(A,B)(\chi u)=P^+F_{1}+P^+F_{2}=(A,B)v_{1}+(A,B)v_{2}.
$$
Здесь функции
\begin{equation}\label{1f33}
v_{1}\in P\bigl(H^{\varphi\varrho^{2q}}(\Omega)\bigr)\quad\mbox{і}\quad v_{2}\in
P\bigl(H^{(s+i+2q)}(\Omega)\bigr)
\end{equation}
являются решениями (единственными) задач
$$
(A,B)v_{1}=P^+F_{1}\in
P^+\bigl(\mathcal{H}^{\varphi}(\Omega,\Gamma)\bigr)\quad\mbox{и}\quad
(A,B)v_{2}=P^+F_{2}\in P^+\bigl(\mathcal{H}^{(s+i)}(\Omega,\Gamma)\bigr).
$$
Теперь из равенства
$$
(A,B)(\chi u)=(A,B)(v_{1}+v_{2})
$$
следует, что
$$
\chi u=v_{1}+(v_{2}+w)\quad\mbox{для деякого}\quad w\in N\subset
C^{\infty}(\overline{\Omega}).
$$
Эта формула с учетом вложений \eqref{1f33} и произвольности функции
$\chi\in\nobreak\Upsilon$ означает истинность заключения импликации \eqref{2f34}.

Таким образом, мы доказали, что эта импликация справедлива для каждого
$i\in\mathbb{N}$. По условию $u\in H^{(s+2q)}(\Omega)$; поэтому посылка импликации
\eqref{2f34} является истинной при $i=1$. Выберем число $p\in\mathbb{N}$ так, чтобы
$s+p>\sigma_{1}(\varphi)$; тогда  в силу \eqref{2f22}
$$
H^{(s+p+2q)}(\Omega)\subset H^{\varphi\varrho^{2q}}(\Omega).
$$
Воспользовавшись импликацией \eqref{2f34} последовательно для значений
$i=1,\,2,\ldots,p$ делаем вывод, что
$$
\chi u\in
H^{\varphi\varrho^{2q}}(\Omega)+H^{(s+p+2q)}(\Omega)=H^{\varphi\varrho^{2q}}(\Omega)
\;\;\mbox{для всех}\;\;\chi\in\Upsilon.
$$
Таким образом, $u\in H^{\varphi\varrho^{2q}}_{\mathrm{loc}}(\Omega_{0},\Gamma_{0})$.

Теорема 4 доказана.

\textbf{\textit{Доказательство теоремы}~5\textit{.}} Произвольно выберем точку
$x\in\Omega_{0}\cup\Gamma_{0}$ и функцию $\chi\in C^{\infty}(\overline{\Omega})$
такую, что $\mathrm{supp}\,\chi\subset\Omega_0\cup\Gamma_0$ и $\chi=1$ в некоторой
окрестности точки~$x$. В~силу теоремы~4, условия \eqref{2f16} и предложения~2, в
котором берем $\alpha(t)\equiv\varphi(t)t^{2q}$, получаем включение
$$
\chi u\in H^{\varphi\rho^{2q}}(\Omega)\subset C^{k}(\overline{\Omega}).
$$
Отсюда с учетом выбора $x$ и $\chi$ следует, что $u\in
C^{k}(\Omega_{0}\cup\Gamma_{0})$. Теорема~5 доказана.

Обоснуем также замечание~1 к этой теореме. Пусть целое $k\geq0$. Предположим, что
для каждого решения $u$, удовлетворяющего условию теоремы~4, выполняется включение
$u\in\nobreak C^{k}(\Omega_{0}\cup\Gamma_{0})$. Покажем, что тогда верно
\eqref{2f16}. Любая функция $u\in H^{\varphi\varrho^{2q}}(\Omega)$ удовлетворяет
условию теоремы~4 и поэтому принадлежит пространству $C^{k}(\Omega_{0})$ согласно
нашего предположения. Отсюда следует вложение
$H^{\varphi\varrho^{2q}}(\Omega_{1})\subset C^{k}(\overline{\Omega}_{1})$, где
$\Omega_{1}$~--- некоторый (произвольно выбранный) открытый шар в $\mathbb{R}^{n}$,
замыкание которого лежит в $\Omega_{0}$. В силу предложения~2 это вложение влечет за
собой \eqref{2f16}.

\textbf{\textit{Доказательство теоремы}~6\textit{.}} Из условий \eqref{2f17}
(включение $f\in H^{\varphi_1}_{\mathrm{loc}}(\Omega,\varnothing)$) и \eqref{2f19}
следует на основании теоремы~5, что $u\in C^{2q}(\Omega)$. Кроме того, из условий
\eqref{2f17} (включение $f\in H^{\varphi_2}(\Omega)$), \eqref{2f18} и \eqref{2f20}
следует в силу той же теоремы, что $u\in C^{m}(\overline{\Omega})$. Теорема~6
доказана.

\end{document}